\documentclass[12pt]{article}
\usepackage{a4}
\usepackage{amsthm}
\usepackage{amsmath}
\usepackage{amssymb}

\newtheorem{thm}{Theorem}
\newtheorem{lem}[thm]{Lemma}

\newtheorem{conj}{Conjecture}

\newcommand{\ind}{\operatorname{ind}}

\begin{document}

\title{A bound on the inducibility of cycles}
\author{Daniel Kr\'al'\thanks{Mathematics Institute, DIMAP and Department of Computer Science, University of Warwick, Coventry CV4 7AL, UK. E-mail: {\tt d.kral@warwick.ac.uk}. The work of this author has received funding from the European Research Council (ERC) under the European Union’s Horizon 2020 research and innovation programme (grant agreement No 648509) and from the Engineering and Physical Sciences Research Council Standard Grant number EP/M025365/1. This publication reflects only its authors' view; the European Research Council Executive Agency is not responsible for any use that may be made of the information it contains.}\and
        Sergey Norin\thanks{Department of Mathematics and Statistics, McGill University, Montreal, Canada. E-mail: {\tt snorin@math.mcgill.ca}. This author was supported by an NSERC grant 418520.}\and
	Jan Volec\thanks{Department of Mathematics and Statistics, McGill University, Montreal, Canada. E-mail: {\tt jan.volec@mcgill.ca}. This author was supported by CRM-ISM fellowship.}}

\date{}

\maketitle

\begin{abstract}
In 1975, Pippenger and Golumbic conjectured that
every $n$-vertex graph has at most $n^k/(k^k-k)$ induced cycles of length $k\ge 5$.
We prove that every $n$-vertex graph has at most $2n^k/k^k$ induced cycles of length $k$.
\end{abstract}

\section{Introduction}
\label{s:intro}

The study of the number of induced copies of a given graph is a classical topic in extremal combinatorics,
which can be traced back to the work of Pippenger and Golumbic~\cite{PipGol75} from 1975.
The \emph{induced density} of a graph $H$ in a graph $G$, which is denoted by $i(H,G)$,
is the number of induced copies of $H$ in $G$ divided by $\binom{|V(G)|}{|V(H)|}$.

A standard averaging argument shows that for all graphs $H$ and $G$ and all integers $|V(H)|\le n<|V(G)|$,
there exists an $n$-vertex graph $G'$ such that $i(H,G')\ge i(H,G)$.
It follows that the sequence $i(H,n)$ is monotone non-increasing in $n$,
and hence it converges for every~$H$.
The \emph{inducibility} of a graph $H$, which is denoted by $\ind(H)$, is the limit of the sequence $i(H,n)$
where $i(H,n)$ is the maximum induced density of $H$ in an $n$-vertex graph.

Pippenger and Golumbic~\cite{PipGol75} showed that
the inducibility of every $k$-vertex graph $H$ is at least $k!/(k^k-k)$ and
conjectured that this bound is tight for a cycle of length $k\ge 5$.
\begin{conj}[Pippenger and Golumbic~\cite{PipGol75}]
\label{conj:main}
The inducibility of a cycle $C_k$ of length $k\ge 5$ is equal to $\frac{k!}{k^k-k}$.
\end{conj}
In the recent years, the flag algebra method of Razborov~\cite{RazFlag}
led to new bounds on the inducibility of small graphs~\cite{BHLP16,Hirst14},
which included the proof of Conjecture~\ref{conj:main} for $k=5$ by Balogh et al.~\cite{BHLP16}.
Other classes of graphs for which the inducibility has been determined
include sufficiently balanced complete multipartite graphs~\cite{BEHJ95,BNT86,BroSid94,PipGol75} and
sufficiently large balanced blow-ups of arbitrary graphs~\cite{HHN14}. 

Motivated by Conjecture~\ref{conj:main},
we study the inducibility of cycles and provide a new upper bound.
In their original paper,
Pippenger and Golumbic~\cite{PipGol75} proved Conjecture~\ref{conj:main} within a multiplicative factor of $2e$,
i.e., they proved that
\[\ind(C_k) \leq \frac{2k!}{k(k-1)^{k-1}} = \left(2e+o(1)\right)\frac{k!}{k^k}.\]
The multiplicative factor $2e$ has recently been improved to $128e/81$ by Hefetz and Tyomkyn~\cite{HefTyo17} and
to $e$ by Pfender and Phillips~\cite{PfePhi}.
Our main result reads as follows.
\begin{thm}
\label{t:main}
Every $n$-vertex graph $G$ contains at most $2n^k/k^k$ induced copies of a cycle $C_k$ of length $k\ge 5$.
\end{thm}
This attains the bound of Conjecture~\ref{conj:main} up to a multiplicative factor of $2$, i.e.,
we show that
\begin{equation}
\ind(C_k) \leq \left(2+o(1)\right)\frac{k!}{k^k}.\label{e:cor}
\end{equation}
We remark that we convinced ourselves that more detailed arguments could be used to improve
the multiplicative factor $2$ in \eqref{e:cor} to $2-\varepsilon$ for some tiny $\varepsilon>0$
but we do not include further details to keep this note short and easily accessible.


\section{Proof of Theorem~\ref{t:main}}

The rest of the paper is devoted to the proof of Theorem~\ref{t:main}.
Fix an $n$-vertex graph $G$ and an integer $k\ge 5$.
Instead of counting the number of induced copies of $C_k$,
we will count the number of $k$-tuples of vertices $(z_1,z_2,z_3,z_4,\ldots,z_k)$ such that
$z_2z_1z_3z_4\cdots z_k$ is an induced cycle of length $k$ in $G$;
we call such a $k$-tuple \emph{good}.
We define a weight $w(D)$ of a good $k$-tuple $D=(z_1,\ldots,z_k)$ as
\[w(D)=\prod_{i=1}^{k}\frac{1}{n_i}\;\mbox{,}\]
where
\begin{itemize}
\item $n_1$ is $n$,
\item $n_2$ is the number of neighbors of $z_1$,
\item $n_3$ is the number of neighbors of $z_1$ that are not neighbors of $z_2$,
\item $n_i$ for $i=4,\ldots,k-1$ is the number of vertices $x$ such that $z_2z_1z_3z_4\cdots z_{i-1}x$ is an induced path of length $i$, and
\item $n_k$ is the number of vertices $x$ such that $z_2z_1z_3z_4\cdots z_{k-1}x$ is an induced cycle of length $k$.
\end{itemize}
In other words, $n_i$ is the number of ways that
we can extend the $(i-1)$-tuple $(z_1,\ldots,z_{i-1})$ by adding a vertex $x$ in a way that
can eventually result in a good $k$-tuple.

The backward induction on $m$ yields that the total weight of good $k$-tuples
starting with the vertices $z_1,\ldots,z_m$ is at most $(n_1\cdots n_m)^{-1}$.
So, we get the following lemma for $m=0$.
We remark that the lemma can also be proven by considering a carefully chosen probability distribution
on some $\ell$-tuples, for $\ell<k$, and good $k$-tuples of vertices of $G$ such that
the probability of choosing a good $k$-tuple $D$ is $w(D)$.

\begin{lem}
\label{l:sum}
The sum of the weights $w(D)$ of all good $k$-tuples $D$ is at most $1$.
\end{lem}

We continue the proof of Theorem~\ref{t:main}.
Consider an induced cycle $v_1v_2v_3\cdots v_k$ of length $k$ in $G$, and
define $D_j$ to be the good $k$-tuple
$(v_j,v_{j-1},v_{j+1},v_{j+2},\ldots,v_{j+k-2})$ for $j=1,\ldots,k$ (indices are modulo $k$).
We will show that
\begin{equation}
\frac{k^k}{4n^k}\le w(D_1)+\cdots+w(D_k)\;\mbox{.}
\label{e:1}
\end{equation}
The inequality \eqref{e:1} implies that
the sum of the $2k$ good $k$-tuples corresponding to a single induced cycle of length $k$ is at least $\frac{k^k}{2n^k}$.
Since the sum of all such $k$-tuples is at most $1$ by Lemma~\ref{l:sum},
the number of induced cycles of length $k$ in $G$ is at most $\frac{2n^k}{k^k}$.
Hence, the proof of Theorem~\ref{t:main} will be completed when we establish \eqref{e:1}.

We now focus on proving \eqref{e:1} and start with applying the AM-GM inequality.
\begin{equation}
\left(\prod_{j=1}^k w(D_j)\right)^{\frac{1}{k}}\le\frac{w(D_1)+\cdots+w(D_k)}{k}
\label{e:AG1}
\end{equation}
Let $n_{j,i}$ be the quantity $n_i$ appearing in the definition of the weight $w(D_j)$.
We obtain the following estimate using the definition of the weight $w(D_j)$, the identity $n_{j,1}=n$ and the AM-GM inequality.
\begin{eqnarray}
\left(\prod_{j=1}^k \frac{1}{w(D_j)}\right)^{\frac{1}{k(k-1)}} 
& = & \left(\prod_{j=1}^k 4n_{j,1}\frac{n_{j,2}}{2}\frac{n_{j,3}}{2}n_{j,4}\cdots n_{j,k}\right)^{\frac{1}{k(k-1)}} \nonumber\\
& = & \left((4n)^k\prod_{j=1}^k \frac{n_{j,2}}{2}\frac{n_{j,3}}{2}n_{j,4}\cdots n_{j,k}\right)^{\frac{1}{k(k-1)}} \nonumber\\
& \le & \frac{(4n)^{\frac{1}{k-1}}}{k(k-1)}\sum_{j=1}^k\frac{n_{j,2}}{2}+\frac{n_{j,3}}{2}+n_{j,4}+\cdots+n_{j,k}\;\mbox{.}
\label{e:AG2}
\end{eqnarray}
We next establish that each vertex $x$ contributes at most $k-1$ to the sum in \eqref{e:AG2}.

We start with showing that each vertex $x$ contributes at most $1$
to the sum $\frac{n_{j,2}}{2}+\frac{n_{j,3}}{2}+n_{j,4}+\cdots+n_{j,k}$ for every $j=1,\ldots,k$.
By symmetry, it is enough to analyze the case $j=1$.
Let $i$ be the smallest index such that $x$ is adjacent to $v_i$.
If $i=1$, then $x$ can contribute only to $n_{1,2}$ and $n_{1,3}$, and if $i=2$, then only to $n_{1,k}$.
If $i=3,\ldots,k-2$, then $x$ can contribute only to $n_{1,i+1}$.
Finally, if $i>k-2$ or $x$ is not adjacent to any vertex $v_i$,
then $x$ does not contribute to any of the summands.
Since the contribution of a vertex $x$ to the sum in \eqref{e:AG2} is at most $1$ for every $j$,
the total contribution of $x$ to the sum in \eqref{e:AG2} is at most $k$;
we improve this bound by $1$ in the next paragraph.

Fix a vertex $x$.
If the vertex $x$ is adjacent to all the vertices $v_1,\ldots,v_k$,
then $x$ contributes $1/2$ to the sum $\frac{n_{j,2}}{2}+\frac{n_{j,3}}{2}+n_{j,4}+\cdots+n_{j,k}$ for every $j$, and
its total contribution to the whole sum in \eqref{e:AG2} is at most $k/2 < k-1$.
Otherwise,
let $i$ be the smallest index such that $x$ is adjacent to $v_{i-1}$ but not to $v_i$ (all indices in this paragraph are modulo $k$).
If $x$ is adjacent to any of the vertices $v_{i+1},\ldots,v_{i+k-4}$ or
it is not adjacent to the vertex $v_{i+k-3}=v_{i-3}$,
then the contribution of $x$ to the sum for $j=i$ is $0$.
Hence, it remains to analyze the following two cases:
\begin{itemize}
\item $x$ is adjacent to the vertices $v_{i-3}$ and $v_{i-1}$ only, and
\item $x$ is adjacent to the vertices $v_{i-3}$, $v_{i-2}$ and $v_{i-1}$ only.
\end{itemize}
In the former case, the contribution of $x$ to the sum for $j=i-2$ is $0$, and
in the latter case, the contribution of $x$ to the sum for $j=i-2$ and for $j=i-1$ is $1/2$.
We conclude that the contribution of each vertex $x$ to the sum in \eqref{e:AG2} is at most $k-1$.

Since the contribution of each vertex $x$ to the sum in \eqref{e:AG2} is at most $k-1$,
the whole sum is at most $n(k-1)$ and we derive the following from \eqref{e:AG2}.
\[
\left(\prod_{j=1}^k \frac{1}{w(D_j)}\right)^{\frac{1}{k(k-1)}} \le
\frac{(4n)^{\frac{1}{k-1}}}{k(k-1)}\cdot n(k-1) = \frac{(4n)^{\frac{1}{k-1}}n}{k}
\]
It follows that
\[
\left(\prod_{j=1}^k \frac{1}{w(D_j)}\right)^{\frac{1}{k}} \le \frac{4n^k}{k^{k-1}}\;\mbox{,}
\]
which is equivalent to
\begin{equation}
\frac{k^{k-1}}{4n^k}\le \left(\prod_{j=1}^k w(D_j)\right)^{\frac{1}{k}}\;\mbox{.}
\label{e:2}
\end{equation}
The desired estimate \eqref{e:1} now follows from \eqref{e:AG1} and \eqref{e:2}.

\section*{Acknowledgments}

The authors would like to thank the anonymous referees for carefully reading the manuscript and for their valuable comments,
which improved the presentation of the results.

\bibliographystyle{bibstyle}
\bibliography{cyclinduce}
\end{document}